\documentclass[11pt]{article}
\usepackage{amssymb,latexsym,amscd,amsmath,amsfonts,enumerate,supertabular}
\usepackage{graphicx}
\usepackage{tabularx}
\usepackage{xcolor}
\usepackage{caption}
\usepackage{epstopdf}
\newtheorem{Theorem}{Theorem}[section] 

\newtheorem{lemma}[Theorem]{Lemma} 
\newtheorem{Proposition}[Theorem]{Proposition}

\numberwithin{equation}{section}
\begin{document}
 \title{{\bf Existence results and numerical method 
 of fourth order convergence
for solving a nonlinear triharmonic equation }}

\author{ Dang Quang A$^{\text a}$,  Nguyen Quoc Hung$^{\text b}$, Vu Vinh Quang $^{\text c}$\\
$^{\text a}$ {\it\small Center for Informatics and Computing, VAST}\\
{\it\small 18 Hoang Quoc Viet, Cau Giay, Hanoi, Vietnam}\\
{\small Email: dangquanga@cic.vast.vn}\\
$^{\text b}$ {\it\small Ha Noi University of Science and Technology,
Ha Noi, Vietnam}\\
{\small Email: hung.nguyenquoc@hust.edu.vn}\\
$^{\text c}$ {\it\small University of Information Technology and Communication,
 Thai Nguyen, Viet Nam}\\
{\small Email: vvquang@ictu.edu.vn}}
\date{ }
%
\maketitle
\begin{abstract}
\small  
{In this work, we consider a boundary value problem for nonlinear triharmonic equation. Due to the  reduction of nonlinear boundary value problems to operator equation for nonlinear terms we establish the existence, uniqueness and positivity of solution. More importantly, we design an iterative method at both continuous and discrete level for numerical solution of the problem. An analysis of actual total error of the obtained discrete solution is made. Some examples demonstrate the applicability of the theoretical results on qualitative aspects and the efficiency of the iterative method.}
\end{abstract}

{\small

}

\section {Introduction}
In this work, we consider the following boundary value problem (BVP) for nonlinear triharmonic equation
\begin{equation}\label{eq1}
\begin{split}
\Delta ^3 u &=f(x,u,\Delta u, \Delta ^2 u), \quad  x \in  \Omega , \\
 u &= 0,\;  \Delta u = 0,\; \Delta ^2 u = 0, \quad  x \in  \Gamma , 
\end{split}
\end{equation}
where $\Omega $ is a  bounded connected domain in $\mathbb{R}^n \; (n \ge 2)$ with the smooth boundary $\Gamma $,  $\Delta $ is the Laplace operator,  $f$ is a continuous function.

The one-dimensional  problem, namely, the problem
\begin{equation}\label{eq3}
\begin{split}
u^{(6)}(x) &=f(x,u(x), u''(x),u^{(4)}(x)), \quad 0 < x < 1,\\
u^{(2i)}(0)&= u^{(2i)}(1)=0 \; (i=0,1,2)
\end{split}
\end{equation}
was studied in many works, e.g., \cite{Al-Ha}, \cite{Bout}, \cite{Khal}, \cite{Lang}, \cite{Noor1}, \cite{Pand}.
The authors of these works were interested only in finding the solution of the problem without attention on the investigation of its qualitative aspects such as the existence, uniqueness and properties of solutions. They either assumed the unique existence of solution or referred  to \cite{Agar} for general results of the existence
and uniqueness of solution for higher order differential equations.
Recently, in \cite{ALong}, Dang and Dang established the existence and uniqueness of solution, and also proposed a  method for the numerical solution.\par
Concerning the triharmonic problem \eqref{eq1}, to the best of our knowledge, the existence of solution  has not been investigated, while the existence of nonlinear biharmonic equation  
\begin{equation}\label{eq1b}
\begin{split}
\Delta^{2}u&=f(x,u,\Delta u), \quad  x \in \Omega,\\
u&=0,\quad \Delta u=0, \quad x \in \Gamma,
\end{split}
\end{equation}
was considered in many works, e.g.,  \cite{An-Liu}, \cite{Choi}, \cite{QA6}, \cite{Hu}, \cite{Liu-Hua}, \cite{Liu-Wang},  \cite{Pao}, \cite{Pei}, \cite{Wang1}. However, for  boundary value problems of sixth order elliptic equations some results of uniqueness of solution are known in  a few works such as \cite{Danet}, \cite{Goyal}, \cite{Schae}. In spite of these facts, there are many works concerning the numerical solution of  nonlinear triharmonic boundary value problems, see e.g., \cite{Ghasemi}, \cite{Mohanty1}-\cite{Mohanty4}. \emph{It should be emphasized that in the above works, under the assumptions that the solution of nonlinear BVPs exist and is unique with sufficient smoothness,  the authors only constructed difference schemes with local truncation error of high order but did not obtain error estimate of the approximate solution. They could obtain error estimate only in the case of linear equation.} \par 

Motivated by the above fact and the importance of sixth order elliptic equation in the modeling of ulcers \cite{Ugail1}, viscous fluid \cite{Lesnic}, geometric  design \cite{Ugail2}, in this paper we  investigate the existence and uniqueness of solution of problem \eqref{eq1} and design an iterative scheme both on continuous and discrete levels for numerical solution of the problem. Especially, we obtain an estimate for total error of the actual approximate solution which is of order four in grid size. The proved theoretical results are validated on examples and the efficiency of the iterative method is shown by numerical experiments.\par
 It is noticed that, the method used here for problem \eqref{eq1} is a further development of our approach to fourth order nonlinear BVPs for ODE and PDE, see, e.g., \cite{QA3}-\cite{QA5}, and for sixth order nonlinear BVPs \cite{ALong}.

\section{Existence and uniqueness of solution}\label{exist}
As in \cite{QA6} for nonlinear biharmonic problem, we shall reduce the problem \eqref{eq1} to 
a fixed point problem.\par
To do this, we define the nonlinear operator $A$ acting in the space of continuous functions $C(\overline{\Omega})$ as follows: for $\varphi (x) \in C(\overline{\Omega})$
\begin{equation}\label{eq2.1}
(A\varphi )(x)=f(x,u(x), \Delta u(x), \Delta ^2 u(x)),
\end{equation}
where $u(x)$ is a solution to the problem
\begin{equation}\label{eq2.2}
\begin{split}
\Delta ^3 u&=\varphi (x), \;  x \in \Omega ,\\
u = 0,  \Delta u &= 0, \Delta ^2 u = 0, \quad  x \in  \Gamma  .
\end{split}
\end{equation}
\begin{Proposition}\label{Prop1}
Function $\varphi (x)$  solves the operator equation
\begin{equation}\label{eq2.3}
A\varphi = \varphi ,
\end{equation}
if and only if the solution $u(x)$ of  problem \eqref{eq2.2} satisfies  problem \eqref{eq1}
\end{Proposition}
{\bf Proof.} The proof of the proposition is similar to that of \cite[Proposition 1]{QA6} (see also \cite[Proposition 2.1]{QA5}).

\noindent{\bf Remark.} The above proposition means that the problem \eqref{eq1} is equivalent to the operator equation \eqref{eq2.3}.\\
 Now we study properties of the operator $A$. 
 First, we note that problem \eqref{eq2.2} can be split into three second order problems
\begin{equation}\label{eq2w}
\begin{aligned}
\Delta w &=\varphi , \quad x \in \Omega , \\ 
w&=0,  \quad x \in \Gamma,
\end{aligned}
\end{equation}
\begin{equation}\label{eq2v}
\begin{aligned}
\Delta v &=w, \quad x \in \Omega , \\ 
v&=0, \quad x \in \Gamma,
\end{aligned}
\end{equation}
\begin{equation}\label{eq2u}
\begin{aligned}
\Delta u&=v, \quad x \in \Omega, \\ 
u&=0, \quad x \in \Gamma.
\end{aligned}
\end{equation} 

\noindent To estimate the solutions of second order problems 
 we need a lemma, which is a corollary of the maximum principle for elliptic equations \cite{Prot} (see also \cite{QA5}, \cite{QA6}). In what follows, we shall use the max-norm 
$$ \| y\|=\max_{x \in \bar \Omega} |y(x)| ,\; \; y(x) \in C(\overline {\Omega} ).$$ 

\begin{lemma}\label{lem1} Let $g(x)$ be a continuous function in $\bar \Omega$. Then for the solution of the problem
\begin{equation}\label{eq5}
\begin{aligned}
-\Delta y &=g, \quad x \in \Omega, \\ 
y &=0,\quad x \in \Gamma ,
\end{aligned}
\end{equation} 
there holds the estimate
\begin{equation}\label{eq6}
\begin{aligned}
\| y\| \leq C_\Omega \| g\|, 
\end{aligned}
\end{equation} 
where $ C_\Omega = R^n/2n$ and $R$ is the radius of the circle which contains the domain $\Omega$. If $\Omega$ is the unit square in $\mathbb{R}^2$ then
\begin{equation}\label{eq7}
\begin{aligned}
\| w\| \leq \frac{1}{8} \| g\|. 
\end{aligned}
\end{equation} 
\end{lemma}

Now we are ready to study  the operator equation \eqref{eq2.2}. Let $M$ be an arbitrary positive number. In the space      $\Omega \times \mathbb{R}^3$  define a bounded domain
\begin{equation}\label{eq8}
\begin{split}
\mathcal {D}_{M}= \Big \{ (x,u,v,w) \mid  x\in \Omega; |u| \leq C_\Omega ^3 M, |v| \leq C_\Omega ^2 M , |w| \leq C_\Omega  M \Big \}
\end{split}
\end{equation}
and  as usual, denote 
\begin{equation*}\label{eq9}
\begin{aligned}
B[O,M]=\{\varphi \in C(\overline {\Omega}) \mid \| \varphi \| \leq M \}. 
\end{aligned}
\end{equation*} 

\begin{Theorem}\label{theorem1}
Assume that there exist numbers $M,L_1, L_2, L_3 \ge 0$ such that 
\begin{enumerate}[(i)]
\item 
\begin{equation}\label{eq12}
|f(x,u,v,v)| \le M \quad  \forall (x,u,v,w) \in \mathcal{D}_M .
\end{equation}
\item 
\begin{equation}\label{eq13}
|f(x,u_2,v_2,w_2)-f(x,v_1,u_1,w_1)| \le L_1 |u_2 -u_1|+L_2 |v_2-v_1|+L_3 |w_2-w_1| \\
\end{equation}
  $ \quad  \forall (x,u_i,v_i,w_i) \in \mathcal{D}_M, \ i=1, 2.$
\item 
\begin{equation}\label{eq14}
q:= (C_\Omega ^2 L_1+ C_\Omega  L_2+ L_3) C_\Omega  < 1.
\end{equation}
\end{enumerate}
Then the BVP \eqref{eq1} possesses a unique solution $u(x) \in C(\bar \Omega)$, satisfying the estimates
\begin{equation}\label{eq15}
\|u\| \le C_\Omega ^3 M,\quad  \|\Delta u\| \le C_\Omega ^2 M, \quad \|\Delta ^2u\| \le C_\Omega  M.
\end{equation}
\end{Theorem}
{\bf Proof.}
Using Proposition \ref{Prop1} and Lemma \ref{lem1} it is easy to prove the above theorem in a similar way as for nonlinear biharmonic problem in \cite[Theorem 1]{QA6} (see also \cite[Theorem 2.3]{QA5}).        $\square$ \par

Now consider a particular case of Theorem \ref{theorem1} concerning the positive solutions.\\
Denote
\begin{equation*}
\begin{split}
\mathcal {D}_{M}^+= \Big \{ (x,u,v,w) \mid  x\in \Omega; 0 \le u \leq C_\Omega ^3 M, - \leq C_\Omega ^2 M \le v \le 0, 0\le w \leq C_\Omega  M \Big \}
\end{split}
\end{equation*}

\begin{Theorem}\label{theorem2}
Assume that there exist numbers $M,L_1, L_2, L_3 \ge 0$ such that 

\begin{equation}\label{eq12a}
-M \le f(x,u,v,v) \le 0 \quad  \forall (x,u,v,w) \in \mathcal{D}_M^+ 
\end{equation}
and the conditions \eqref{eq13}, \eqref{eq14} are satisfied in $\mathcal {D}_{M}^+$.
Then the BVP \eqref{eq1} possesses a unique positive solution $u(x) \in C(\bar \Omega)$, satisfying the estimates
\begin{equation}\label{eq15a}
0 \le u \le C_\Omega ^3 M, \quad - C_\Omega ^2 M \le \Delta u \le 0,\quad  0\le \Delta ^2u \le C_\Omega  M.
\end{equation}
\end{Theorem}

\noindent {\bf Proof.} 
The proof of the theorem is similar to that of Theorem \ref{theorem1} if replace the ball $B[0,M]$ by the strip
\begin{equation*}
S_M= \{ \varphi\in C(\bar \Omega) | \; -M \le \varphi(x) \le 0 \}.
\end{equation*}

\section{Iterative method at continuous level}
Consider the following iterative method for solving
the problem \eqref{eq1}  based on the successive approximation of the fixed point of the operator $A$:\\
\noindent 
i) Given an initial  approximation $\varphi_0 \in B[O,M]$, for example, 
\begin{equation}\label{eq3.1}
\varphi _0(x)=f(x,0,0,0), \quad x \in \Omega.
\end{equation}
ii) Knowing $\varphi _k$ $(k=0,1,2,...)$ solve sequentially three second order problems
\begin{equation}\label{eq3.2}
\begin{aligned}
\Delta w_k &=\varphi_k, \quad x \in \Omega , \\ 
w_k&=0,\quad x \in \Gamma,
\end{aligned}
\end{equation}
\begin{equation}\label{eq3.3}
\begin{aligned}
\Delta v_k &=w_k, \quad x \in \Omega , \\ 
v_k&=0,\quad x \in \Gamma,
\end{aligned}
\end{equation}
\begin{equation}\label{eq3.4}
\begin{aligned}
\Delta u_k &=v_k, \quad x \in \Omega , \\ 
u_k&=0,\quad x \in \Gamma .
\end{aligned}
\end{equation}
iii) Calculate the new approximation
\begin{equation}\label{eq3.5}
\varphi_{k+1}(x)= f(x,u_k(x),v_k(x),w_k(x)).
\end{equation}
\begin{Theorem}\label{theo3}
Assume that all the assumptions of Theorem \ref{theorem1} are satisfied. Then the iterative method \eqref{eq3.1}-\eqref{eq3.5} converges geometrically and there holds the error estimate 
\begin{equation}\label{eq3.6}
\| u_k-u\| \leq \frac{C_\Omega^2 q^k}{1-q}\| \varphi _1 -\varphi _0\|,
\end{equation}
$u$ being the exact solution of the problem \eqref{eq1} and $q$ being given by \eqref{eq14}.
\end{Theorem} 
\noindent {Proof.}
Indeed, as the iterative method \eqref{eq3.1}-\eqref{eq3.5} is the successive approximation method for 
the fixed point of the operator $A$ with the initial approximation  $\varphi _0 \in B[O,M]$, we have
\begin{equation}\label{eq3.7}
||\varphi _k - \varphi || \leq \frac{q^k}{1-q}||\varphi _1 - \varphi _0 ||.
\end{equation}
Next, from the estimate
$$ \|u_k-u\| \le C_\Omega ^3 \|\varphi _k -\varphi \|
$$
which follows from Lemma \ref{lem1} applied consecutively to the problems

\begin{equation*}
\begin{aligned}
\Delta (u_k-u) &=v_k-v, \quad x \in \Omega , \\ 
u_k-u&=0,\quad x \in \Gamma ,
\end{aligned}
\end{equation*}
\begin{equation*}
\begin{aligned}
\Delta (v_k -v)&=w_k -w, \quad x \in \Omega , \\ 
v_k&=0,\quad x \in \Gamma,
\end{aligned}
\end{equation*}
\begin{equation*}
\begin{aligned}
\Delta (w_k-w) &=\varphi_k -\varphi, \quad x \in \Omega , \\ 
w_k&=0,\quad x \in \Gamma,
\end{aligned}
\end{equation*}
and \eqref{eq3.7} we obtain \eqref{eq3.6}. Thus, the theorem is proved. $\square$

\section{Iterative method at discrete level}

We restrict ourselves to two-dimensional problems, and consider the problem in the unit square
$\bar{\Omega}=[0, 1] \times [0, 1]$. We cover $\overline{\Omega}$ by the uniform grid 
 
$$\overline{\omega}_h = \Big \{(x_1,x_2)|\; x_1=ih_1, x_2=jh_2, i=\overline{0,m}, j = \overline{0, l} \Big \},$$ 
where $h_1=1/m, h_2=1/l.$ Denote by $\Omega _h$ and $\gamma_h$ the set of interior points and the set of boundary points of $\overline{\omega}_h$, respectively.\\
For solving the Poisson problems \eqref{eq3.2}-\eqref{eq3.4} at each iterative step we shall use finite difference schemes of fourth order of accuracy. For this purpose, denote by $\Phi_k(x), W_k(x),V_k(x), U_k(x)$ the grid functions defined on the grid $\overline{\omega}$ and approximating the functions $\varphi_k(x) , w_k(x) , v_k(x) , u_k(x) $ on this grid.\par
Before describing iterative method we need a result from the theory of difference schemes \cite{Sam1}. Consider the Dirichlet problem for Poisson equation
\begin{equation}\label{4eq1}
\begin{split}
\Delta y &=\psi (x), \; x\in \Omega,\\
y&=\mu (x), x\in \Gamma .
\end{split}
\end{equation}
Denote by $Y(x)$ a grid function approximating the function $y(x)$ on the grid $\overline{\omega} _h$ and recall the following notations: 
\begin{equation*}\label{4eq2}
\begin{aligned}
\Lambda _1 Y&=\dfrac{Y_{i-1,j}-2Y_{ij}+Y_{i+1,j}}{h_1^2},\;
\Lambda _2 Y=\dfrac{Y_{i,j-1}-2Y_{ij}+Y_{i,j+1}}{h_2^2},\\
\Lambda Y &=(\Lambda _1 + \Lambda _2)Y,\quad 
\Lambda ^* Y = \Lambda Y + \dfrac{h_1^2+h_2^2}{12}\Lambda _1\Lambda _2 Y, \\
\psi^* &= \psi + \dfrac{h_1^2}{12}\Lambda _1\psi + \dfrac{h_2^2}{12}\Lambda _2\psi ,
\end{aligned}
\end{equation*}
where $Y_{ij}=Y(ih_1,jh_2).$
\begin{lemma}\label{lemma2} (See \cite[Sec. 4.5]{Sam1})
For the boundary value problem \eqref{4eq1} the difference scheme
\begin{equation}\label{4eq3}
\begin{split}
\Lambda ^* Y &= \psi^*,\;  x\in \omega _h\\
Y&=\mu (x), \; x\in \gamma _h 
\end{split}
\end{equation}
is of fourth order accuracy, when $y \in C^{(6)}(\Omega)$, i.e.
$\|U-u\|_h =O(h^4)$,
where
\begin{equation}\label{4eq4}
\|U-u\|_h =\max _{x\in \overline{\omega}_h}|Y(x)-y(x)|.
\end{equation}
If replace $\psi ^*$ by $\hat{\psi} = \psi ^* +O(h^4)$ then the order of accuracy is not changed.
\end{lemma}

Now consider the following discrete iterative method:
\begin{enumerate}
\item Given 
\begin{equation}\label{4eq5}
\Phi_0(x) =f(x,0,0,0),\; x \in \omega_h.
\end{equation}
\item Knowing $\Phi_k$ in $\omega_h$ $(k=0,1,...)$ solve consecutively three difference problems
\begin{equation}\label{4eq6}
\begin{split}
\Lambda^* W_k&={\Phi_k}^*,\quad x \in \omega_h,\\
{W_k}&=0,\quad x \in \gamma_h,
\end{split}
\end{equation}
\begin{equation}\label{4eq7}
\begin{split}
\Lambda^* V_k&={W_k}^*,\quad x \in \omega_h,\\
{V_k}&=0,\quad x \in \gamma_h,
\end{split}
\end{equation}
\begin{equation}\label{4eq8}
\begin{split}
\Lambda^* U_k&={V_k}^*,\quad x \in \omega_h,\\
U_k&=0,\quad x \in \gamma_h,
\end{split}
\end{equation}
\item Compute the new approximation
\begin{equation}\label{4eq9}
\Phi_{k+1}=f(x,U_k,V_k,W_k),\quad x \in \omega_h.
\end{equation}
\end {enumerate}

\begin{Proposition}\label{Prop2}
Assume that the function $f(x,u,v,w)$ is continuous and has all continuous partial derivatives up to sixth order in the domain $\mathcal{D}_M$. Then for the functions $u_k(x), v_k(x),w_k(x)$ constructed by the iterative method \eqref{eq3.1}- \eqref{eq3.5} there hold
$$ w_k(x) \in C^{(6)}(\Omega), \; v_k(x) \in C^{(8)}(\Omega), \; u_k(x) \in C^{(10)}(\Omega).
$$
\end{Proposition}
\noindent {Proof.} The proposition is easily proved by induction. $\square $

\begin{Proposition}\label{Prop3}
Under the assumptions of Proposition \ref{Prop2},  for  $k=0,1,2,...$  we have the estimates
\begin{align}
\|\Phi_k -\varphi _k \|_h &=O(h^4), \label{4eq10} \\
\|U_k -u _k \|_h &= O(h^4),  \label{4eq11} \\
\|V_k -v _k \|_h &= O(h^4),  \label{4eq12} \\
\|W_k -w _k \|_h &= O(h^4).  \label{4eq13}
\end{align}
\end{Proposition}
\noindent {Proof.} The proposition will be proved by induction.\\
When $k=0$ the difference scheme for the problem
\begin{equation}\label{4eq14}
\begin{aligned}
\Delta w_0 &=\varphi_0, \quad x \in \Omega, \\ 
w_0 &=0,\quad x \in \Gamma ,
\end{aligned}
\end{equation}
is 
\begin{equation}\label{4eq15}
\begin{split}
\Lambda^* W_0&={\Phi _0}^*,\quad x \in \omega_h,\\
{W_0}&=0,\quad x \in \gamma_h,
\end{split}
\end{equation}
where
$${\Phi_0 }^* = \Phi_0 + \dfrac{h_1^2}{12}\Lambda _1\Phi_0 + \dfrac{h_2^2}{12}\Lambda _2\Phi_0, \quad x\in \omega_h . $$
As $\Phi_0 =\varphi_0 =f(x,0,0,0)$ we have ${\Phi_0 }^* = {\varphi }^* $. By Lemma \ref{lemma2} we obtain the estimate
\begin{equation}\label{4eq16}
\|W_0 -w _0 \|_h = O(h^4). 
\end{equation}
Next, consider the difference scheme
\begin{equation}\label{4eq17}
\begin{split}
\Lambda^* V_0&={W _0}^*,\quad x \in \omega_h,\\
{V_0}&=0,\quad x \in \gamma_h
\end{split}
\end{equation}
for the problem
\begin{equation}\label{4eq18}
\begin{aligned}
\Delta v_0 &=w_0, \quad x \in \Omega, \\ 
v_0 &=0,\quad x \in \Gamma ,
\end{aligned}
\end{equation}
where
$${W_0 }^* = W_0 + \dfrac{h_1^2}{12}\Lambda _1 W_0 + \dfrac{h_2^2}{12}\Lambda _2 W_0.
$$
From \eqref{4eq16} it follows 
$$ {W_0 }^* = {w_0 }^* +O(h^4).$$
Once again, by Lemma \ref{lemma2} we obtain
$$ \|V_0 -v _0 \|_h = O(h^4).  $$
Analogously, we have
$$ \|U_0 -u _0 \|_h = O(h^4).  $$
Now, suppose that the estimates \eqref{4eq10}-\eqref{4eq13} are valid for $k \ge 0$. We must prove them for $k+1$. Since $\varphi _{k+1}$ and $\Phi _{k+1}$ are calculated by \eqref{eq3.5} and \eqref{4eq9}, respectively, and the function 
$f(x,u,v,w)$ satisfies the Lipschitz condition in the variables $u,v$ and $w$, we have
$$ |\Phi _{k+1}-\varphi_{k+1}| \le L_1 |U_k  -u_k|+L_2 |V_k  -v_k| +L_3 |W_k  -w_k|.$$
From the above estimate and the hypothesis of induction it follows
$$ \|\Phi_{k+1} -\varphi_{k+1} \|_h = O(h^4).  $$
Now,  making the same argument as above for $k=0$, we consecutively obtain the estimates
\begin{align*}
\|W_{k+1} -w _{k+1} \|_h &= O(h^4),    \\
\|\Phi_{k+1} -\varphi _{k+1} \|_h &=O(h^4),  \\
\|V_{k+1} -v _{k+1} \|_h &= O(h^4),   \\
\|U_{k+1} -u _{k+1} \|_h &= O(h^4).
\end{align*}
Thus, the proposition is proved. $\square $

\begin{Theorem}\label{theo4}
Under the assumptions of Theorem 1, for the solution of the problem  \eqref{eq1} obtained by the discrete iterative method \eqref{4eq5}-\eqref{4eq9},  there hold the estimates
\begin{align}
\|U_k -u  \|_h &= C_\Omega ^3  p_kd +O(h^4),  \label{4eq21} \\
\|V_k -v  \|_h &=C_\Omega ^2  p_kd + O(h^4),  \label{4eq22} \\
\|W_k -w  \|_h &= C_\Omega   p_kd + O(h^4),  \label{4eq23}
\end{align}
where 
\begin{equation}\label{4eq24}
p_k=\frac{q^k}{1-q}, \; d=\| \varphi _1 - \varphi _0 \|,\; v=\Delta u ,\; w=\Delta ^2u.
\end{equation}
\end{Theorem}
{\bf Proof.} Representing 
$$U_k -u =(U_k-u_k )+ (u_k-u) $$
we have
$$\|U_k -u \|_h\le \|U_k-u_k \|_h+ \|u_k-u\|_h .$$
In view of the estimates \eqref{eq3.6} and \eqref{4eq11} we obtain \eqref{4eq21}.\\
Analogously, we have \eqref{4eq22} and \eqref{4eq23}. The theorem is proved. $\square$

\section{Examples}
To demonstrate the applicability of the existence results in Section \ref{exist} and show the efficiency of the iterative method in the previous section we shall consider several examples. As was said in the beginning of the previous section, all examples will be considered in the computational domain $\Omega =[0, 1] \times [0,1]$ with the boundary $\Gamma$.
For testing the convergence of the proposed iterative method we perform some experiments for the cases, where  exact solutions are known and for the cases where exact solutions are not known.
In all examples,  the iterative process \eqref{4eq5}-\eqref{4eq9} will carried out until
\begin{equation}\label{5eq1}
\|\Phi_{k} -\Phi_{k-1} \|_h \le TOL ,
\end{equation}
where $TOL$ is given accuracy. For solving the discrete problems \eqref{4eq6}-\eqref{4eq8} we use the cyclic reduction method \cite{Sam2}, which is one of the direct methods for grid equations.
\\

{\bf{Example 1.} } Consider the problem
\begin{equation*}
\begin{split}
\Delta ^3 u &=- \pi^4 \sin (\pi x_1) \sin(\pi x_2) + (\Delta u)^2 -u \Delta ^2 u +\sin(\Delta ^2u-\Delta ^2u^*), \quad x\in \Omega , \\
 u &=  \Delta u =\Delta ^2 u = 0, \quad x \in \Gamma , 
\end{split}
\end{equation*}
where 
$$u^*= \frac{1}{8\pi ^2}\sin (\pi x_1) \sin(\pi x_2). $$
For this example
\begin{equation*}
f(x,u,v,w)= - \pi ^4 \sin (\pi x_1) \sin(\pi x_2) +v^2 -u w +\sin(w-\Delta ^2u^*).
\end{equation*}
It is possible to verify that for $M=104$ there holds $|f(x,u,v,w)| \le M$ in 
\begin{equation*}
\begin{split}
\mathcal {D}_{M}= \Big \{ (x,u,v,w) \mid  x\in \Omega; |u| \leq  M/8^3, |v| \leq  M/8^2 , |w| \leq  M/8 \Big \}.
\end{split}
\end{equation*}
Besides, in $\mathcal{D}_M$, all other assumptions of Theorem \ref{theorem1} are also satisfied with 
$L_1=M/8=13,\ L_2=2M/8^2= 3.25, \ L_3= M/8+1=1.2031 $ and 
$q=0.2266 <1$.
Therefore, by this theorem the problem has a unique solution satisfying the estimate $|u(x)|\le M/8^3=0.2031$.
It is easy verify that this unique solution is the above function $u^*(x_1,x_2)$.

For testing the convergence of the discrete iterative method we perform numerical experiments on computer LENOVO, 64-bit Operating System (Win 10), Intel Core I5, 1.8 GHz, 8 GB RAM  with stopping criterion 
\eqref{5eq1} for different $TOL$.
The results of computation are reported in Tables \ref{Tab1} and \ref{Tab2}. 

\begin{table}[!ht]
\centering
\setlength{\tabcolsep}{0.6cm}
 \caption[smallcaption]{The results of computation of Example 1 for $TOL=10^{-6}$ }
\label{Tab1}
\begin{tabular}{ ccccc} 
\hline
$\text{Grid }$ & $K$ & $ E(K)  $ & $ Order  $ & $Time$ (in sec.)\\
\hline
16$\times$ 16     & 6 & 1.5399e-07   &3.9979   &0.0781 \\
32$\times$ 32     & 6 & 9.6387e-09   &3.8682   & 0.1094 \\
64$\times$ 64     & 6 & 6.0006e-10   & 4.1042  & 0.3438\\
128$\times$ 128   & 6 & 3.4891e-11   &6.2757   & 4.1094\\
{\bf 256$\times$ 256 }  & 6 & 4.5035e-13   & -2.5581 & 15.8125 \\
512$\times$ 512   & 6 & 2.6523e-12   & 0.0364    & 54.5469\\
1024$\times$ 1024 & 6 & 2.5862e-12   &  0.2803      & 273.7969\\
2048$\times$ 2048 & 6 & 2.1296e-12   &      & 2.6920e+03     \\
\hline
\end{tabular}
\end{table}
\noindent From the two first columns of Table \ref{Tab1} we see that the number of iterations $K$ for achieving  the same tolerance $TOL= 10^{-6}$ 
is independent of the grid size, and it is $K=6$. The two next columns have the following meaning: 
$E^h(K)= \|U^h_K-u^* \|_h$, $Order$ is the order of convergence calculated by the formula
$$ Order=\log _2 \frac{\|U^h_K-u^*\|_h}{\|U^{h/2}_K-u^*\|_{h/2}}.
$$
In the above formula the superscripts $h$ and $h/2$ of $U$ mean that $U$ is computed on the grid with the corresponding grid sizes. \\
From Table \ref{Tab1} we see that, first, when the number of grid points increases from $16 \times 16$ to $256 \times 256$ the accuracy $E(K)$ increases with the rate close to $4$ or more. After that the accuracy almost remains unchanged. This phenomenon completely agrees with the estimate \eqref{4eq21} of Theorem \ref{theo4} because after the two terms of errors are  the same then further decrease of the second term does not improve the total error. This phenomenon also will be observed in Table \ref{Tab2}. The results of computation for other tolerances also support this assertion.

\begin{table}[!ht]
\centering
\setlength{\tabcolsep}{0.6cm}
 \caption[smallcaption]{The results of computation of Example 1 for $TOL=10^{-8}$ }
\label{Tab2}
\begin{tabular}{ ccccc} 
\hline
$\text{Grid }$ & $K$ & $ E(K)  $ & $ Order  $ & $Time$ (in sec.)\\
\hline
16$\times$ 16     & 8 & 1.5399e-07   &3.9975   &0.0781 \\
32$\times$ 32     & 8 &  9.6414e-09   &3.9994   & 0.1719 \\
64$\times$ 64     & 8 & 6.0284e-10   & 4.0003  & 0.5156\\
128$\times$ 128   & 8 & 3.7670e-11   &4.0161   & 5.1719\\
256$\times$ 256   & 8 & 2.3282e-12   & 4.2041 & 21.6250 \\
{\bf 512$\times$ 512 }  & 8 & 1.2632e-13   & -0.6279    & 72.6250\\
1024$\times$ 1024 & 8 & 1.9521e-13   &  -1.7359      & 368.4844\\
2048$\times$ 2048 & 8 & 6.5020e-13   &      & 3.5572e+03
     \\
\hline
\end{tabular}
\end{table}

\noindent {\bf{Example 2.} }Consider the problem
\begin{equation*}
\begin{split}
\Delta ^3 u &=- \sin (\pi x_1) \sin(\pi x_2) + u (\Delta u) -\dfrac{1}{2} \Delta ^2 u , \quad x\in \Omega , \\
 u &=  \Delta u =\Delta ^2 u = 0, \quad x \in \Gamma . 
\end{split}
\end{equation*}
In this example
\begin{equation*}
f(x,u,v,w)= - \pi \sin (\pi x_1) \sin(\pi x_2) +uv -\dfrac{1}{2} w.
\end{equation*}
Choosing $M=1.5$, it is easy to verify that in the domain
\begin{equation*} 
\begin{split}
\mathcal {D}_{M}^+= \Big \{ (x,u,v,w) \mid  x\in \Omega; 0 \le u \leq \dfrac{1}{8^3} M, 
-\dfrac{1}{8^2}M \le v \le 0, 0\le w \le \dfrac{1}{8}M \Big \}
\end{split} 
\end{equation*}
all the conditions of Theorem \ref{theorem2} are satisfied .
Consequently, the problem has a unique positive solution. \\
The results of convergence for $TOL=10^{-8}$ are given  in Table \ref{Tab3}.
\begin{table}[!ht]
\centering
\setlength{\tabcolsep}{0.6cm}
 \caption[smallcaption]{The results of computation of Example 2 for $TOL=10^{-8}$}
\label{Tab3}
\begin{tabular}{ cccc} 
\hline
$\text{Grid }$ & $K$ & $ e^h(K)  $ & $ Order  $\\
\hline
16$\times$ 16   & 7 & 3.4346 e-14 &3.9974  \\
32$\times$ 32   & 7 & 3.4348 e-14 &3.9993\\
64$\times$ 64   & 7 & 3.4348 e-14 &3.9993\\
128$\times$ 128 & 7 & 3.4348 e-14 & 4.0043\\
256$\times$ 256 & 7 & 3.4348 e-14 & 5.0082\\
{\bf 512$\times$ 512} & 7 & 3.4348 e-14 &-2.8568\\
1024$\times$ 1024 & 7 & 3.4348 e-14   &     \\
2048$\times$ 2048 & 7 & 3.4348 e-14   &     \\
\hline
\end{tabular}
\end{table}
Here, in the case when the exact solution is unknown, the deviation between two successive iterations $ e^h(K)$ and $Order$ of convergence are calculated by the formulas
\begin{align*}
 e^h(K)&= \|U^h_{K}-U^h_{K-1} \|_h,\\
 Order &=\log _2 \frac{\|U^h_K-U^{h/2}_K\|_h}{\|U^{h/2}_K-U^{h/4}_K\|_{h/2}}.
\end{align*}

\begin{figure}
\begin{center}
\includegraphics[height=6cm,width=10cm]{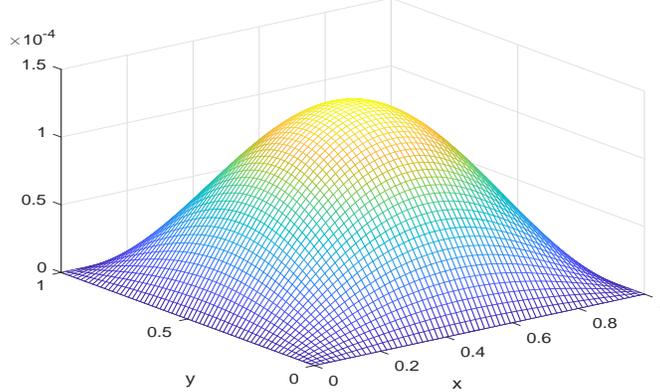}
\caption{The graph of the approximate solution in Example $2$. }
\label{Fig1}
\end{center}
\end{figure}

The graph of the approximate solution computed on the grid $64 \times 64 $ is depicted in Figure \ref{Fig1}.\\

Remark that  Examples 1 and 2 demonstrate the validity of the theoretical results of the existence, uniqueness and positivity of solution of the problem \eqref{eq1} and the convergence of the proposed iterative method \eqref{4eq5}-\eqref{4eq9}. Below we consider an example, where the conditions of Theorem \ref{theorem1} are not satisfied. In this case we are not sure about the existence and uniqueness of solution, and the convergence of the iterative method. However, the results of computation show that the iterative method converges to a solution. This means that Theorem \ref{theorem1} only gives sufficient conditions for existence, uniqueness of solution and convergence of the iterative method.\\

\noindent {\bf{Example 3.} }Consider the problem
\begin{equation*}
\begin{split}
\Delta ^3 u &=-100+50u^2-(\Delta u )^3+ \Delta ^2u , \quad x\in \Omega , \\
 u &=  \Delta u =\Delta ^2 u = 0, \quad x \in \Gamma . 
\end{split}
\end{equation*}
In this example
\begin{equation*}
f(x,u,v,w)= -100+50u^2-v^3+w^2.
\end{equation*}
 We cannot choose $M$ such that the first assumption of Theorem \ref{theorem1} is satisfied. Indeed, the inequality
$$100+50\Big(\dfrac{M}{8^3}\Big)^2+\Big(\dfrac{M}{8^2}\Big)^3+\Big(\dfrac{M}{8}\Big)^2 \le M
$$ 
 has no solution. Therefore, the existence of solution of the problem is not guaranteed. However, the iterative method \eqref{4eq5}-\eqref{4eq9} converges. The results of computation are given in Table 4 and the graph of computed approximation on grid $64 \times 64$ is depicted in Figure 2.

\begin{table}[!ht]
\centering
\setlength{\tabcolsep}{0.6cm}
 \caption[smallcaption]{The results of computation of Example 3 for $TOL=10^{-8}$}
\label{Tab4}
\begin{tabular}{ cccc} 
\hline
$\text{Grid }$ & $K$ & $ e^h(K)  $ & $ Order  $\\
\hline
16$\times$ 16   & 33 & 6.2588 e-13 &3.9981  \\
32$\times$ 32   & 33 & 6.2515 e-13 &3.9996\\
64$\times$ 64   & 33 & 6.2511 e-13 &3.9998\\
128$\times$ 128 & 33 & 6.2508 e-13 & 4.0010\\
256$\times$ 256 & 33 & 6.2509 e-13 & 4.5615\\
{\bf 512$\times$ 512} & 33 & 6.2514 e-13 &-0.2625\\
1024$\times$ 1024 & 33 & 6.2515 e-13   &     \\
2048$\times$ 2048 & 33 & 6.2518 e-13    &     \\
\hline
\end{tabular}
\end{table}

\begin{figure}
\begin{center}
\includegraphics[height=6cm,width=10cm]{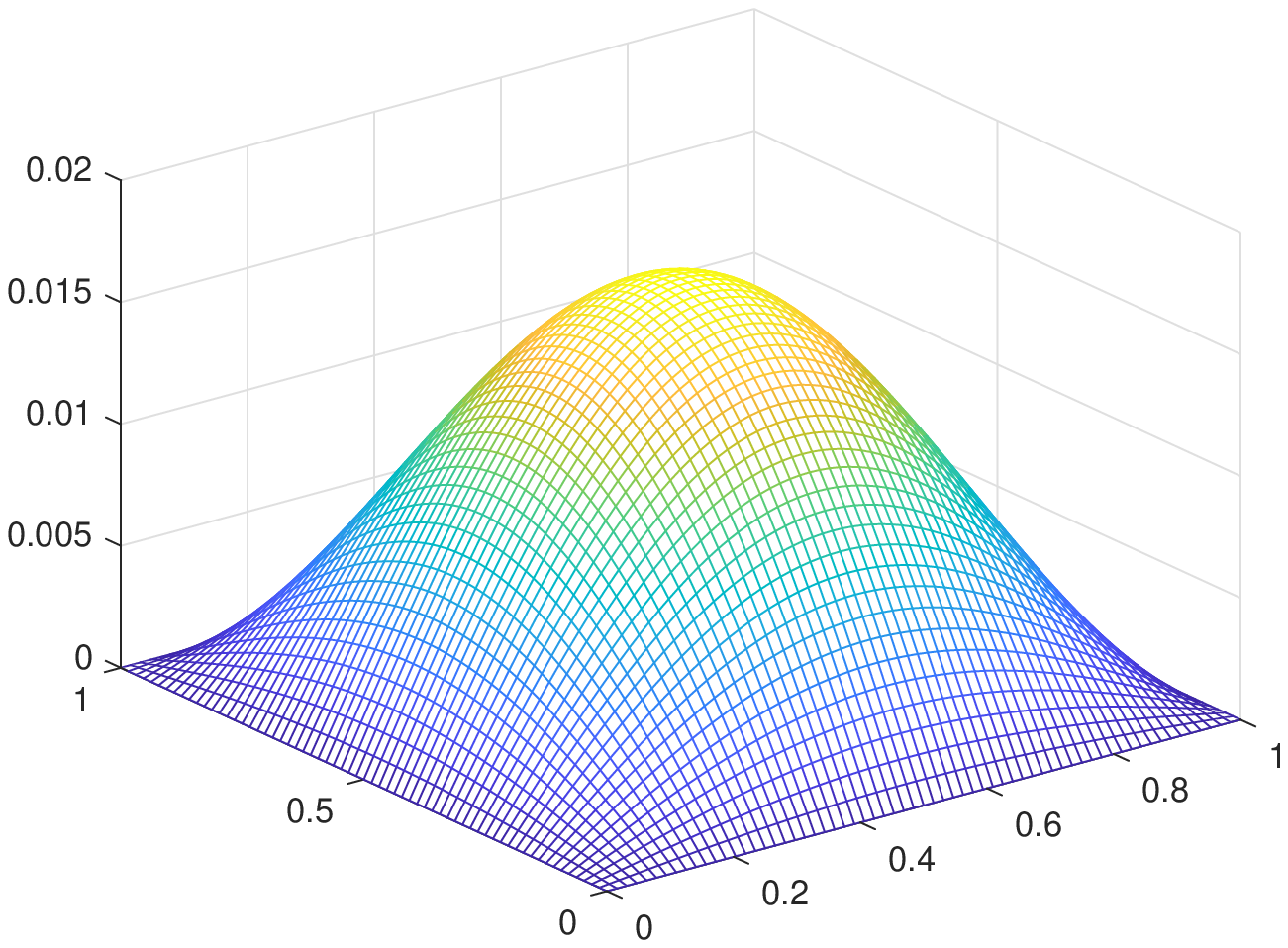}
\caption{The graph of the approximate solution in Example $3$. }
\label{FigExam3}
\end{center}
\end{figure}

\noindent {\bf{Example 4.} }Consider the problem with nonhomogeneous boundary conditions
\begin{equation*}
\begin{split}
\Delta ^3 u &=f(x,u,\Delta u, \Delta ^2 u), \quad  x \in  \Omega , \\
 u &= g_1,\;  \Delta u = g_2,\; \Delta ^2 u = g_3, \quad  x \in  \Gamma . 
\end{split}
\end{equation*}
In this case, to solve the problem numerically we use the following discrete iterative method:
\begin{enumerate}
\item Given 
\begin{equation}\label{eqIT1}
\Phi_0(x) =f(x,0,0,0),\; x \in \omega_h.
\end{equation}
\item Knowing $\Phi_k$ in $\omega_h$ $(k=0,1,...)$ solve consecutively three second order difference problems
\begin{equation}\label{eqIT2}
\begin{split}
\Lambda^* W_k&={\Phi_k}^*,\quad x \in \omega_h,\\
{W_k}&=g_3,\quad x \in \gamma_h,
\end{split}
\end{equation}
\begin{equation}\label{eqIT3}
\begin{split}
\Lambda^* V_k&={W_k}^*,\quad x \in \omega_h,\\
{V_k}&=g_2,\quad x \in \gamma_h,
\end{split}
\end{equation}
\begin{equation}\label{eqIT4}
\begin{split}
\Lambda^* U_k&={V_k}^*,\quad x \in \omega_h,\\
U_k&=g_1.\quad x \in \gamma_h,
\end{split}
\end{equation}
\item Compute the new approximation
\begin{equation}\label{eqIT5}
\Phi_{k+1}=f(x,U_k,V_k,W_k),\quad x \in \omega_h.
\end{equation}
\end {enumerate}
This discrete iterative method is expected to be convergent of fourth order, too. 
Below we give an numerical example illustrating the fourth order convergence of the above iterative method.\\
Consider the equation
\begin{equation*}
\Delta ^3u =4e^{x_1+x_2}+\sin(e^{x_1+x_2}-u)-\cos (2e^{x_1+x_2}-\Delta u)+ \Delta ^2 u +1
\end{equation*}
with the exact solution $u = e^{x_1+x_2}$. The boundary conditions are calculated from this exact solution. The results of computation  are reported in Tables \ref{Tab5} and \ref{Tab6}.

\begin{table}[!ht]
\centering
\setlength{\tabcolsep}{0.6cm}
 \caption[smallcaption]{The results of computation of Example 4 for $TOL=10^{-6}$}
\label{Tab5}
\begin{tabular}{ cccc} 
\hline
$\text{Grid }$ & $K$ & $ E(K)  $ & $ Order  $\\
\hline
16$\times$ 16   & 8 & 1.6538e-08 & 4.0004     \\
32$\times$ 32   & 8 & 1.0333e-09 &4.0345    \\
64$\times$ 64   & 8 & 6.3054e-11 &4.3822    \\
128$\times$ 128 & 8 & 3.0238e-12 & 1.0156   \\
256$\times$ 256 & 8 & 1.4957e-12 & -1.7774   \\
 512$\times$ 512 & 8 & 5.1275e-12 &-2.7983   \\
1024$\times$ 1024 & 8 & 3.5667e-11   &     \\
\hline
\end{tabular}
\end{table}

\begin{table}[!ht]
\centering
\setlength{\tabcolsep}{0.6cm}
 \caption[smallcaption]{The results of computation of Example 4 for $TOL=10^{-8}$}
\label{Tab6}
\begin{tabular}{ cccc} 
\hline
$\text{Grid }$ & $K$ & $ E(K)  $ & $ Order  $\\
\hline
16$\times$ 16   & 10 & 1.6539e-08 & 3.9983     \\
32$\times$ 32   & 10 & 1.0349e-09 & 3.9984    \\
64$\times$ 64   & 10 & 6.4758e-11 & 3.7714   \\
128$\times$ 128 & 10 & 4.7424e-12 & 0.8431   \\
256$\times$ 256 & 10 & 2.6437e-12 & -0.5815    \\
 512$\times$ 512 & 10 & 3.9559e-12 &  -3.0981  \\
1024$\times$ 1024 & 10 & 3.3875e-11   &     \\
\hline
\end{tabular}
\end{table}

\begin{figure}
\begin{center}
\includegraphics[height=6cm,width=10cm]{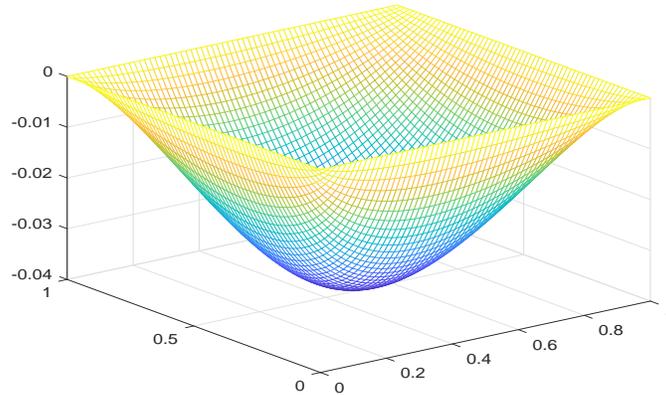}
\caption{The graph of the approximate solution in Example $5$. }
\label{Fig3}
\end{center}
\end{figure}

\noindent {\bf{Example 5.} } Now consider an example of using the iterative method \eqref{eqIT1}-\eqref{eqIT5} for a problem with nonhomogeneous boundary conditions, for which the exact solution is not known. The problem is
\begin{equation*}
\begin{split}
\Delta ^3 u &=8 e^{-(x_1^2 +x_2^2)}+\sin(u)+\dfrac{1}{1+(\Delta  u)^2}+e^{-\Delta ^2 u}, \quad  x \in  \Omega , \\
 u &= \sin (x_1(1-x_1)x_2(1-x_2)),\;  \Delta u = e^{-(x_1^2 +x_2^2)},\; \Delta ^2 u = \cos(x_1 x_2), \quad  x \in  \Gamma . 
\end{split}
\end{equation*}
The result of computation on the grid $64 \times 64$ with $\|\Phi_{k} -\Phi_{k-1} \|_h \le 10^{-6}$ is that the iterative process stops after $K=7$ iterations and $\|U_K- U_{K-1}\|_h=4.4222e-11$. The graph of the approximate solution is depicted in Figure \ref{Fig3}.

\section{Conclusion}
In this work, by reducing the original boundary value problem of nonlinear triharmonic equation to an operator equation for the nonlinear term we have established the existence, uniqueness and positivity of solution of it. And more importantly, we have designed a numerical iterative method which is proved to be of fourth order of accuracy in total. This is achieved due to the use of difference schemes of the same order for the Poisson equation at each iteration. 
%
Several examples, where the exact solutions of the problem are known or are not known, demonstrate the validity of obtained theoretical results and the efficiency of the proposed iterative method.  \par
In the future we shall use the technique of this paper in combination with the boundary operator method in \cite{QA2002} to consider the nonlinear triharmonic equation with other boundary conditions. This is a perspective direction of research.

\newpage



\begin{thebibliography}{00}

\bibitem{Agar} R. P. Agarwal, Boundary value problems for higher order differential equations, World Scientific, Singapore, 1986.

\bibitem{Al-Ha}W. Al-Hayani, Adomian decomposition method with Green’s function for sixth-order
boundary value problems, Computers and Mathematics with Applications 61 (2011) 1567-1575.

\bibitem {An-Liu} Y. An, R. Liu, Existence of nontrivial solutions of an asymptotically linear fourth-order elliptic equation, Nonlinear Analysis, 68 (2008) 3325-3331.

\bibitem{Bout}A. Boutayeb, E. H. Twizell, Numerical methods for the solution of special sixth order boundary-value problems, International Journal of Computer Mathematics, 45(3-4)(1992), 207-223.

\bibitem {Choi} Q.H. Choi, T. Jung, A fourth order nonlinear elliptic equation with jumping nonlinearity, Houston J. Math., 24 (1998) 735-756.

\bibitem {QA2002} Q. A Dang, Iterative Method for Solving a Boundary Value Problem for Triharnonic Equation, Vietnam Journal of Mathematics, 30:l (2002) 71-78.

\bibitem {ALong}	Q. A. Dang, Q. L. Dang, A simple efficient method for solving sixth-order nonlinear boundary value problems, Comp. Appl. Math. (2018) 37 (Suppl 1): 16. 

%

\bibitem{QA3} Q. A. Dang, T. K. Q. Ngo, “Existence results and iterative method for solving the caltilever beam equation with fully nonlinear term”, Nonlinear Analysis: Real World Applications,  36 (2017) 56-68.

\bibitem{QA4} Q. A. Dang, Q. L. Dang, T. K. Q. Ngo, A novel efficient method for nonlinear boundary value problems, Numerical Algorithms, 76 (2017), 427-439.

\bibitem{QA5} Q. A. Dang, T. H. Nguyen, “Existence result and iterative method for solving a nonlinear biharmonic equation of Kirchhoff type”, Computers \& Mathematics with Applications, 76 (2018), 11-22.

\bibitem{QA6} Q. A. Dang, H. H. Truong, T. H. Nguyen, T. K. Q. Ngo, Solving a nonlinear biharmonic boundary value problem, Journal of Computer Science and Cybernetics, 33 (4) (2017), 308-324.

\bibitem {Danet} C. P. Danet, Uniqueness in some higher order elliptic boundary value problems in $n$-dimensional domains, Electronic Journal of Qualitative Theory of Differential Equations, 2011, No. 54, 1-12.

\bibitem {Ghasemi} M. Ghasemi, On the numerical solution of high order multi-dimensional elliptic PDEs, Computers and Mathematics with Applications 76 (2018) 1228-1245.

\bibitem {Goyal} S. Goyal and V. Goyal, Liouville-Type and Uniqueness Results for a Class of Sixth-Order Elliptic Equations, Journal of Mathematical Analysis and Applications, 139 (1989), 586-599 .

\bibitem {Hu} S. Hu, L. Wang, Existence of nontrivial solutions for fourth-order asymptically linear elliptic equations, Nonlinear Analysis 94 (2014) 120-132.

\bibitem{Khal}M.Khalid, M. Sultana and F. Zaidi,Numerical Solution of Sixth-Order Differential
Equations Arising in Astrophysics by Neural Network, International Journal of Computer Applications, 107 (2014): 6, 1-6.

\bibitem{Lang} F. G. Lang and X. P. Xu, An Effective Method for Numerical Solution and Numerical 
Derivatives for Sixth Order Two-Point Boundary Value Problems, Computational Mathematics and Mathematical Physics, 2015, Vol. 55, No. 5, pp. 811-822.

\bibitem{Lesnic} D. Lesnic, On the boundary integral equations for a two-dimensional slowly rotating highly viscous fluid flow, Adv. Appl. Math. Mech. 1 (2009) 140-150.

\bibitem {Liu-Hua} X. Liu, Y. Huang , On sign-changing solution for a fourth-order asymptotically linear elliptic problem, Nonlinear Analysis, 72 (2010) 2271-2276.

\bibitem {Liu-Wang} Y. Liu, Z.P. Wang, Biharmonic equations with asymptotically linear nonlinearities, Acta Math. Sci. 27B (3) (2007) 549-560.

\bibitem {Mohanty1} R. K. Mohanty, M K Jain and B N Mishra, A compact discretization of $O(h^4)$ for two-dimensional nonlinear triharmonic equations, Physica Scripta 84 (2011) ID: 025002, doi:10.1088/0031-8949/84/02/025002 
  
\bibitem {Mohanty2} R.K. Mohanty, Single cell compact finite difference discretizations of order two and four for multidimensional triharmonic problems, Numer. Meth. Partial Diff. Eq., 26 (2010) 1420-1426

\bibitem {Mohanty3} R. K. Mohanty, M. K. Jain and B. N. Mishra, A Novel Numerical Method of $O(h^4)$ for Three-Dimensional Non-Linear Triharmonic Equations, Commun. Comput. Phys., 12 (5) (2012) 1417-1433.

\bibitem {Mohanty4} B. N. Mishra and M. K. Mohanty, Single Cell Numerov Type Discretization for 2D Biharmonic and Triharmonic Equations on Uniqual Mesh, Journal of Mathematical and Computational Science, 3 (2013) 242-253.

\bibitem{Noor1} M. A. Noor, K. I. Noor, S. T. Mohyud-Din, Variational iteration method for solving sixth-order boundary value problems, Commun Nonlinear Sci Numer Simulat 14 (2009) 2571-2580.

\bibitem{Pand} P. K. Pandey, Fourth order finite difference method for sixth order boundary value problems, Comput. Math.Math. Phys. 53 (2013) 57-62.

\bibitem{Pao} C.V. Pao, On fourth-order elliptic boundary value problems, Proc. Amer. Math. Soc. 128 (2000) 1023-1030.

\bibitem {Pei} R. Pei, Multiple solutions for biharmonic equations with asymptotically linear nonlinearities,   Boundary Value Problems, (2010), V. 2010, Article ID 241518.

\bibitem {Prot} M. H. Protter, H. F. Weinberger, Maximum Principles in Differential Equations, Springer, 1984.
\bibitem {Sam1} A. A. Samarskii,  The Theory of Difference Schemes, New York, Marcel Dekker, 2001.

\bibitem {Sam2} A. A. Samarskii, E. Nikolaev, Numerical methods for grid equation, Vol. 1, Direct methods, Birkhauser, Basel, 1989.

\bibitem {Schae} Schaefer, P.W. Uniqueness in some higher order elliptic boundary value problems. Journal of Applied Mathematics and Physics (ZAMP) 29, 693–697 (1978). https://doi.org/10.1007/BF01601494

\bibitem {Ugail1} H. Ugail, M. Wilson, Modelling of oedemus limbs and venous ulcers using partial differential equations, Theor. Biol. Med. Model. 2 (2005) 1-28.

\bibitem {Ugail2} H. Ugail, Partial Differential Equations for Geometric Design, Springer, 2011.

\bibitem {Wang1} F. Wang, Y. An, Existence and multiplicity of solutions for a fourth-order elliptic equation,  Boundary Value Problems, (2012), 2012:6.

\end{thebibliography}
\end{document}